\documentclass[a4paper,12pt]{amsart}
\usepackage{amssymb}
\usepackage{yfonts}

\def\hbar{\bar{h}}

\def\iso{\buildrel \sim\over\to}

\def\CD{{\mathcal{D}}}

\def\CI{{\mathcal{I}}}

\def\CO{{\mathcal{O}}}

\def\BQ{{\mathbf{Q}}}

\def\BZ{{\mathbf{Z}}}

\def\can{{\mathrm{can}}}

\def\Spec{\operatorname{Spec}\nolimits}

\newtheorem{thm}{Theorem}[section]
\newtheorem{lemma}[thm]{Lemma}
\newtheorem{cor}[thm]{Corollary}
\newtheorem{prop}[thm]{Proposition}

\newtheorem{conj}[thm]{Conjecture}
\theoremstyle{definition}
\newtheorem{rem}[thm]{Remark}
\newtheorem{example}[thm]{Example}

\numberwithin{equation}{section}

\usepackage{comment}
\usepackage{stmaryrd}
\usepackage{imakeidx}

\input xypic
\xyoption{all}

\makeindex
\makeindex[name=ter,title=Terms,columns=2]

\includecomment{comment}

\title{Cactus groups and Lusztig's asymptotic algebra}

\author{Rapha\"el Rouquier and Noah White}
\address{Rapha\"el Rouquier~: UCLA Mathematics Department, Los Angeles, CA 90095-1555, USA.}
\email{rouquier@math.ucla.edu}
\address{Noah White~: Mathematical Sciences Institute, Australian National University, Canberra, ACT,
Australia}
\email{noah.white@anu.edu.au}
\thanks{The first author gratefully acknowledges support from the NSF
(grant DMS-1702305) and the
Simons Foundation (grant \#376202).}

\begin{document}
\maketitle
\begin{abstract}
	We construct a morphism from the cactus group associated with
a Coxeter group to the group of invertible elements of
Lusztig's asymptotic algebra. This relates to
the cactus group action on elements of Coxeter groups defined in \cite{Lo,Bo2}
and we propose a conjecture on how to fully recover those actions.
\end{abstract}

\section{Introduction}

Cactus groups are "crystal limits" of braid groups, originally introduced
implicitely in type $A$ by Drinfeld \cite{Dr} and explicitely in \cite{HeKa}.
The braid group action on tensor powers of representations of quantum groups
becomes an action of the cactus group on tensor powers of a crystal.

Cactus groups are (orbifold) fundamental groups of the
Deligne-Mumford compactification of the moduli space of genus $0$ real
curves with marked points \cite{De,DaJaSc}.

\smallskip
Cactus groups have been generalized to other Coxeter groups. They can
be defined by generators and relations and, for finite Coxeter groups,
are orbifold fundamental
groups of real points of the wonderful compactification of projectivized
hyerplane complements \cite{DaJaSc}.

\smallskip
In \cite{Lo}, Losev constructed an action of the cactus group associated
to a Weyl group on the set of elements of the Weyl group in terms of the
combinatorics of perverse self-equivalences of the category $\CO$ of
a complex semi-simple Lie algebra.
Bonnaf\'e \cite{Bo2}
generalized this construction to all Coxeter
groups and unequal parameters, with a direct algebraic approach using
the Hecke algebra.
The cactus group orbits are proven to be unions of Kazhdan-Lusztig cells
(actual Kazhdan-Lusztig cells in type $A$).

Cactus groups are expected to play a role in the geometry
of ramification of Calogero-Moser spaces, which conjecturally provides
another construction of Kazhdan-Lusztig cells \cite{BoRou1,BoRou2}.

\smallskip
In this article, following a suggestion of Etingof, we
start from Drinfeld's unitarization trick, providing a morphism
from the cactus group to the (completed) braid group. This enables us
to obtain a direct connection from the cactus group
to the Hecke algebra and to Lusztig's asymptotic algebra. Our work follows
Etingof's proposal and answers some of his questions.

A key point in our approach is a characterization of the element of
the Hecke algebra that corresponds to the longest element of a finite
Coxeter group, given a suitable
isomorphism of the Hecke algebra with the group algebra of the Coxeter group
(\S 3).

\medskip
We thank Pavel Etingof very warmly for introducing us to cactus groups and
for sharing his insights.  We also thank C\'edric Bonnaf\'e for
useful discussions.

\section{Braid groups and cactus groups}
\subsection{Coxeter groups}
Let $(W,S)$ be a Coxeter group with $S$ finite. Let $(m_{st})_{s,t\in S}$ be the associated
Coxeter matrix.

Given $I$ a subset of $S$, let $W_I$ be the subgroup
of $W$ generated by $I$. When $W_I$ is finite, we say that $I$ is spherical
and we denote by $w_I$ the longest element of $W_I$.
We put $w_0=w_S$ when $W$ is finite.

\subsection{Braid groups}

We denote by $B_W$ the braid group of $W$. It is generated by $(\beta_s)_{s\in S}$ and
relations
$$\underbrace{\beta_s\beta_t\beta_s\cdots}_{m_{st}\text{ terms }}=
\underbrace{\beta_t\beta_s\beta_t\cdots}_{m_{st}\text{ terms }}.$$

\smallskip
There is a surjective morphism of groups
$$p_W:B_W\to W,\ \beta_s\mapsto s.$$
Its kernel $P_W$ is the pure braid group.

\smallskip
Let $w\in W$. Given $w=s_1\cdots s_n$ a reduced decomposition, the element
$\beta_w=\beta_{s_1}\cdots \beta_{s_n}$ of $B_W$ is independent of the reduced decomposition.

\smallskip

Given $I$ a subset of $S$ such that $W_I$ is finite, we put
$\beta_I=\beta_{w_I}$. Note that $\beta_I^2\in P_W$.

\subsection{Completions}
Let $\CI$ be the augmentation ideal of the group algebra $\BQ[P_W]$ (the kernel of the 
algebra morphism $\BQ[P_W]\to\BQ,\ P_W\ni g\mapsto 1$). Let
$\widehat{\BQ[P_W]}$ be the completion of $\BQ[P_W]$ at $\CI$. This is a complete cocommutative
Hopf algebra and we denote by $\widehat{\BQ[P_W]}^*$ its topological dual.
The prounipotent completion of $P_W$ is $\hat{P}_W=\Spec(\widehat{\BQ[P_W]}^*)$.
Note that given $g\in P_W$ and $\alpha\in\BQ$, we have
an element $g^\alpha\in\hat{P}_W$ corresponding to
$\sum_{n\ge 0}{\alpha\choose n}(g-1)^n\in \widehat{\BQ[P_W]}$.

\smallskip
Let $\CI'$ be the kernel of $\BQ[p_W]:\BQ[B_W]\to\BQ[W],\ B_W\ni g\mapsto p_W(g)$.
This is the two-sided ideal of $\BQ[B_W]$ generated by $\CI$. We denote by
$\widehat{\BQ[B_W]}$ be the completion of $\BQ[B_W]$ at $\CI'$ and we put
$\hat{B}_W=\Spec(\widehat{\BQ[B_W]}^*)$. This is a proalgebraic group, the connected
component of the identity is $\hat{P}_W$ and $p_W$ extends to a surjective morphism
of groups $\hat{B}_W\to W$ with kernel $\hat{P}_W$.

\subsection{Cactus group}

\smallskip
The cactus group $C_W$ is the group generated by $(\gamma_I)_{\substack{
	I\subset S \text{ spherical}}}$ with
relations
$$\gamma_I^2=1$$
$$\gamma_I\gamma_J=\gamma_{I\cup J} \text{ if }W_{I\cup J}=W_I\times W_J$$
$$\gamma_I\gamma_J=\gamma_J\gamma_{w_J(I)} \text{ if }I\subset J.$$

\smallskip
There is a surjective morphism of groups
$$\pi_W:C_W\to W,\ \gamma_I\mapsto w_I.$$

\smallskip
Note that $C_W$ is generated by those elements $\gamma_I$ such that 
$(W_I,I)$ is a finite irreducible Coxeter group.

\subsection{Cactus to completed braids}

In type $A$,
the following result is in \cite[proof of Theorem 3.14]{EtHeKaRa}, based on
Drinfeld's unitarization trick \cite{Dr}.

\begin{prop}
	\label{pr:morphism}
	The assignment
	$\gamma_I\mapsto \beta_I (\beta_I^2)^{-1/2}$ for $I\subset S$ spherical
	defines a morphism of
	groups $\phi:C_W\to \hat{B}_W$. We have $p_W\circ\phi=\pi_W$.
\end{prop}

\begin{proof}
	Since $\beta_I$ commutes with $\beta_I^2$, it commutes with
	$(\beta_I^2)^{-1/2}$, hence $(\beta_I (\beta_I^2)^{-1/2})^2=1$.

	Similarly, if $W_{I\cup J}=W_I\times W_I$, then $\beta_I$ and $\beta_J$ commute,
	hence $\beta_I (\beta_I^2)^{-1/2}$ and $\beta_J (\beta_J^2)^{-1/2}$ commute.

	Finally, when $I\subset J$, we have $\beta_J^{-1}\beta_I\beta_J=\beta_{w_J(I)}$,
	hence $\beta_J^{-1}(\beta_I^2)^{-1/2}\beta_J=(\beta_{w_J(I)}^2)^{-1/2}$.
	It follows that
	$\beta_I (\beta_I^2)^{-1/2}\cdot \beta_J (\beta_J^2)^{-1/2}=
	\beta_J (\beta_J^2)^{-1/2}\cdot \beta_{w_J(I)} (\beta_{w_J(I)}^2)^{-1/2}$ since
	$(\beta_J^2)^{-1/2}$ commutes with $\beta_{w_J(I)}$ and with 
	$(\beta_{w_J(I)}^2)^{-1/2}$.
\end{proof}

\begin{rem}
	The map of \cite[proof of Theorem 3.14]{EtHeKaRa} and \cite{Dr}
	is defined using a different set of generators for $C_W$. That
	new set of generators can be defined for an arbitrary $W$ as follows.

	Let $F$ be the set of pairs $(I,s)$ where $I\subset S$, $s\in I$ and
	where $W_I$ is finite and irreducible (i.e. its Coxeter diagram
	is connected). Let $(I,s)\in F$. We put
	$\gamma_{I,s}=\gamma_I\gamma_{I\setminus\{s\}}$,
	$\beta_{I,s}'=\beta_I\beta_{I\setminus \{s\}}^{-1}$ and
	$\beta_{I,s}''=\beta_{I\setminus \{s\}}^{-1}\beta_I$,
	where $\gamma_\emptyset=1$ and $\beta_\emptyset=1$. 
	The set $\{\gamma_{I,s}\}_{(I,s)\in F}$
	generates $C_W$ and we have
	$\phi(\gamma_{I,s})=\beta_{I,s}'(\beta_{I,s}''\beta_{I,s}')^{-1/2}$.
\end{rem}

\section{Hecke algebra}

Let $H_W$ be the Hecke algebra of $W$. This is the quotient of the group ring
$\BZ[v^{\pm 1}]B_W$ by the ideal generated by
$(\beta_s-v)(\beta_s+v^{-1})$ for $s\in S$. We denote by
$\kappa:\BZ[v^{\pm 1}]B_W\to H_W$ the quotient map and we put
$T_s=\kappa(\beta_s)$. Note that the 
composition
$$\BZ[v^{\pm 1}]B_W\xrightarrow{\kappa} H_W\xrightarrow{v\to 1}\BZ W$$
is given by $\beta\mapsto p_W(\beta)$.

\medskip
Let $R$ be the completion of $\BQ[v^{\pm 1}]$ at the ideal $(v-1)$.
We have $\kappa(\CI)\subset (v-1)\BQ[v^{\pm 1}] H_W$. It follows that $\kappa$ induces a morphism
of proalgebraic groups (still denoted by $\kappa$) $\hat{B}_W\to (RH_W)^\times$
and we have a commutative diagram
\begin{equation}
	\label{eq:commute}
	\xymatrix{
	C_W \ar[d]_{\pi_W} \ar[r]^-{\phi} &\hat{B}_W \ar[r]^-{\kappa} & (RH_W)^\times
	\ar[d]^{v\to 1} \\
	W \ar@{^{(}->}[rr]_-{\can} && (\BQ W)^\times
}
\end{equation}

Given $I\subset S$, we have a corresponding commutative diagram
with $W$ replaced by $W_I$. It is a subdiagram of the commutative diagram (\ref{eq:commute}).

\begin{lemma}
	\label{le:rank1}
	If $S=\{s\}$, then $\kappa\circ\phi(\gamma_s)=\frac{1-v^2}{1+v^2}+\frac{2v}{1+v^2}T_s$.
\end{lemma}

\begin{proof}
	Let $a=\kappa\circ\phi(\gamma_s)$. We have $a_{|v=1}=s$ and $a^2=1$. This shows that
	$a=\frac{1-v^2}{1+v^2}+\frac{2v}{1+v^2}T_s$.
\end{proof}

\begin{prop}
	\label{pr:defw_I}
Assume $W$ is finite.
	There exists a unique element $\tilde{w}_0$ of $RH_W$ such that 
	\begin{itemize}
		\item $\tilde{w}_0^2=1$
		\item $(\tilde{w}_0)_{|v=1}=w_0$
		\item $\tilde{w}_0T_s\tilde{w}_0^{-1}=T_{w_0sw_0}$ for all
			$s\in S$.
	\end{itemize}
\end{prop}

\begin{proof}
		Given $s\in S$, we have $\gamma_S\gamma_{s}\gamma_S=
		\gamma_{w_0sw_0}$, hence
	$$\bigl(\kappa\circ\phi(\gamma_S)\bigr)T_s\bigl(\kappa\circ\phi(\gamma_S)\bigr)^{-1}=
	T_{w_0sw_0}$$	
	by Lemma \ref{le:rank1}.
	The commutativity of the diagram (\ref{eq:commute}) shows that
	$\tilde{w}_0=\kappa\circ\phi(\gamma_S)$ satisfies the
	properties of the proposition.

	\smallskip
	Consider now an element $h\in RH_W$ satisfying the properties
	of $\tilde{w}_0$.
	Let $z=h\tilde{w}_0\in Z(RH_W)$. We have $z^2=1$  and
	$z_{|v=1}=1$. There is an isomorphism of $R$-algebras
	$i:Z(RH_W)\iso R^n$ for some $n$. Since $i(z)^2=1$, it
	follows that $i(z)\in \BQ^n$. In addition, $i(z)_{|z=1}=1$, hence
	$i(z)=1$. So, we have shown that $z=1$.
	This proves the unicity of the element $\tilde{w}_0$.
\end{proof}

Given $I\subset S$ spherical, we denote by $\tilde{w}_I$ the 
element of $RH_I\subset RH_W$ defined for $W_I$ as in Proposition
\ref{pr:defw_I}. The proof of Proposition \ref{pr:defw_I} shows the following.

\begin{prop}
	\label{pr:charinHecke}
	Given $I\subset S$ with $W_I$ finite, we have
	$\kappa\circ\phi(\gamma_I)=\tilde{w}_I$.
\end{prop}

The next proposition follows from Proposition \ref{pr:charinHecke} and
Lemma \ref{le:rank1}.

\begin{prop}
	\label{pr:cases}
	Given $s\in S$, we have $\tilde{w}_s=\frac{1-v^2}{1+v^2}+\frac{2v}{1+v^2}T_s$.
\end{prop}

The next result is immediate.
\begin{prop}
	Assume $W$ is finite.
Let $\lambda:RW\iso RH_W$ be an isomorphism
of $R$-algebras such that
\begin{itemize}
	\item $\lambda_{|v=1}$ is the identity
	\item $\lambda(w_0sw_0)=T_{w_0}\lambda(s)T_{w_0}^{-1}$ for all $s\in S$.
\end{itemize}

	We have $\lambda(w_0)=\tilde{w}_0$.
\end{prop}

Note that Lusztig has constructed an explicit isomorphism with these
properties which is already defined over $\BQ[v^{\pm 1}]_{(v-1)}$ 
\cite[Theorem 3.1]{Lu1}. As a consequence,
we have the following result answering a question of Etingof.

\begin{cor}
	\label{co:field}
Given $I\subset S$ with $W_I$ finite, we
	have $\tilde{w}_I\in \BQ[v^{\pm 1}]_{(v-1)}H_W$.
\end{cor}

\section{Asymptotic algebra}

Let $h\mapsto \bar{h}$ the $\BZ$-algebra involution of $H_W$ given by
$\bar{v}=v^{-1}$ and $\overline{T_s}=T_s^{-1}$. There is a unique family
$(C_w)_{w\in W}$ of elements of $H_W$ such that $\overline{C_w}=C_w$ and
$C_w-T_w\in\bigoplus_{w'<w}\BZ[v^{\pm 1}]T_{w'}$. This is the
Kazhdan-Lusztig basis of $H_W$. We refer to \cite{Bo1} for basics of
Kazhdan-Lusztig and Lusztig theory.

\smallskip
Given $w,w',w''\in W$, we define $h_{w,w',w''}\in\BZ[v^{\pm 1}]$ so that
$C_wC_{w'}=\sum_{w''\in W}h_{w,w',w''}C_{w''}$.

We put $a(w)=-\min_{w',w''\in W}\deg(h_{w',w'',w})$ and we put
$\gamma_{w,w',w''}=(v^{a(w'')}h_{w,w',(w'')^{-1}})_{|v=0}$.

Lusztig's asymptotic $J$-ring is the $\BZ$-algebra with basis
$\{t_w\}_{w\in W}$ and multiplication given by
$t_wt_{w'}=\sum_{w''\in W}\gamma_{w,w',w''}t_{w''}$.
There is a subset $\CD$ of $W$ such that $1=\sum_{d\in\CD}t_d$.

\medskip

There is a morphism of $\BZ[v^{\pm 1}]$-algebras
$$\psi:H_W\to \BZ[v^{\pm 1}]J_W,\
C_w\mapsto \sum_{y\in W,d\in\CD, a(d)=a(y)} h_{w,d,y}t_y.$$

By Corollary \ref{co:field}, the map $\kappa\circ \phi$ takes values in
$\BQ[v^{\pm 1}]_{(v-1)}$ and we put
$f=\psi\circ\kappa\circ \phi:C_W\to (\BQ[v^{\pm 1}]J_W)^\times$.

\medskip
We have a commutative diagram where all the maps are equivariant for the
automorphisms of $(W,S)$.

$$
	\xymatrix{
		C_W \ar[ddrr]_f \ar@{..>}[drr] \ar[r]^-{\phi}  &\hat{B}_W \ar[r]^-{\kappa} & (RH_W)^\times
	 \\
	 && (\BQ[v^{\pm 1}]_{(v-1)}H_W)^\times \ar@{^{(}->}[u] \ar[d]^\psi\\
	 && (\BQ[v^{\pm 1}]_{(v-1)}J_W)^\times \\
}
$$

\medskip
The following result describing the image of $\gamma_S\in C_W$ in
the asymptotic algebra is our main result.

\begin{thm}
	\label{th:main}
	Assume $W$ is finite.
	We have $f(\gamma_S)=\sum_{d\in\CD}(-1)^{\ell(w_0)+a(w_0d)}
	t_{w_0 d}$.
\end{thm}

\begin{proof}
	Let $t=\sum_{d\in\CD}(-1)^{\ell(w_0)+a(w_0d)}t_{w_0 d}$. We have
	$t^2=1$ and $tt_wt=t_{w_0ww_0}$ \cite[\S 2.9]{Lu2} (cf also
	\cite[Example 19.3.3]{Bo1}), hence
	$\psi^{-1}(t)^2=1$ and $\psi^{-1}(t)T_w\psi^{-1}(t)=
	T_{w_0ww_0}$ for all $w\in W$
	since $\psi$ is equivariant for the diagram
	automorphism induced by $w_0$.

	Since $\psi^{-1}(t)-T_{w_0}\in (v-1)H_W$ 
	\cite[Corollary 2.8]{Lu2}, we have
	$\psi^{-1}(t)-\tilde{w}_0\in (v-1)H_W$, hence
	$\psi^{-1}(t)=\tilde{w}_0$ by Proposition \ref{pr:defw_I}.

	Proposition \ref{pr:charinHecke} completes the proof of the Theorem.
\end{proof}

It follows from \cite[Theorem 3.1]{Ma} (cf \cite[Example 19.3.3]{Bo1})
that given $w\in W$, there is 
$\sigma_W(w)\in W$ such that 
\begin{equation}
\label{eq:Matthas}
f(\gamma_S)t_w=(-1)^{a(w_0w)} t_{\sigma_W(w)} \text{ and }
t_w f(\gamma_S)=(-1)^{a(w_0w)} t_{w_0\sigma_W(w)w_0}.
	\end{equation}

The map $\sigma_W$ is an involution of the set $W$.

\medskip

In general, the image of $f(\gamma_I)$ in $J_W$ is difficult
to describe explictly. When $I=\{s\}$, this is equivalent to
the description of $\psi(C_s)$ (cf Proposition \ref{pr:cases}). We propose
that $f(\gamma_I)$ has an explicit description modulo $v\BZ[v]J_W$ and
that left and right multiplication by $f(\gamma_I)$ give the cactus
group actions on $W$ of Losev and Bonnaf\'e.

\begin{conj}
	\label{co:conj}
	\begin{itemize}
		\item[(1)]
	Given $I\subset S$ spherical, we have
			$f(\gamma_I)\in\BZ[v]_{(v(v-1))}J_W$.
\item[(2)] Denote by $\bar{f}:C_W\to J_W^\times$ the composition of $f$ with
	specialization of $v$ to $0$. Given $I\subset S$ spherical,
			we have $$\bar{f}(\gamma_I)=\sum_{d\in\CD}
			(-1)^{a(w_Id_1)} t_{\sigma_{W_I}(d_1)d_2}$$
			where $d=d_1d_2$ with $d_1\in W_I$ and $d_2$
			of minimal length in $W_Iw$.
\item[(3)] Consider $I\subset S$ spherical and $w\in W$.
Write $w=xy$ with $x\in W_I$ and $y$ of minimal length in $wW_I$ (resp.
$w=yx$ with $x\in W_I$ and $y$ of minimal length in $W_Iw$). We have
		$$\bar{f}(\gamma_I)t_w=(-1)^{a(w_Ix)} t_{\sigma_{W_I}(x)y}
			\text{ and }
			t_w\bar{f}(\gamma_I)=(-1)^{a(w_Ix)} t_{yw_I\sigma_{W_I}(x)w_I}.$$
	\end{itemize}
\end{conj}

\begin{prop}
	Conjecture \ref{co:conj} holds when $I=S$ or $|I|=1$.
\end{prop}

\begin{proof}
	Assume $I=\{s\}$. We have $\tilde{w}_s=-1+\frac{2v}{1+v^2}C_s$, hence
$$f(\gamma_S)t_w=-t_w+\frac{2v}{1+v^2}
	\sum_{w'\in W,a(w')=a(w)} h_{s,w,w'}t_{w'}.$$
	We have 
	$$C_sC_w=
	\begin{cases}
		(v+v^{-1})C_w & \text{ if }sw<w \\
		\sum_{w'\in W, sw'<w'}\mu_{s,w,w'} C_{w'}& \text { otherwise,}
	\end{cases}$$
	where $\mu_{s,w,w'}=h_{s,w,w'}\in\BZ$ when $sw>w$.
	It follows that
	$$f(\gamma_S)t_w=\begin{cases}
		t_w & \text{ if }sw<w \\
		-t_w+\frac{2v}{1+v^2}\sum_{w'\in W, a(w')=a(w), sw'<w'}
		\mu_{s,w,w'} t_{w'}& \text{ otherwise.}
	\end{cases}$$
	This shows that
	$f(\gamma_s)=\sum_{d\in\CD}f(\gamma_S)t_d\in \BZ[v]_{(v(v-1))}J_W$
	and this also shows the first statement of (3). The second statement of (3) is proven
	similarly by considering $t_wf(\gamma_S)$.
	Also we obtain $\bar{f}(\gamma_s)=\sum_{d\in\CD}(-1)^{\delta_{sd<d}}
	t_d$, which shows (2).

	\medskip
	Assume now $I=S$. Statements (1) and (2) are given by Theorem
	\ref{th:main} while (3) is statement (\ref{eq:Matthas}).
\end{proof}

$$
	\xymatrix{
		C_W \ar[ddrr]_f \ar@{..>}[drr]\ar[r]^-{\phi} \ar@/_2pc/ @{-->} [dddrr]^?
		\ar@/_2pc/ @{-->} [ddddrr]_{\bar{f}}
		&\hat{B}_W \ar[r]^-{\kappa} & (RH_W)^\times
	 \\
	 && (\BQ[v^{\pm 1}]_{(v-1)}H_W)^\times \ar@{^{(}->} [u] \ar[d]^\psi\\
	 && (\BQ[v^{\pm 1}]_{(v-1)}J_W)^\times \\
	 && (\BZ[v]_{(v(v-1))}J_W)^\times \ar@{^{(}->}[u]\ar[d]_{v\to 0} \\
	 && (J_W)^\times 
}
$$

\begin{rem}
	We do not know if the map $f:C_W\to (\BQ[v^{\pm 1}]_{(v-1)}J_W)^\times$ is faithful or
	not, and the faithfulness of $\phi:C_W\to\hat{B}_W$ does not seem
	to be known either.
	On the other hand, cactus groups are known to be linear when 
	$W$ is finite \cite{Yu}.
\end{rem}

\begin{example}
Fix $m\ge 3$ and consider the dihedral Coxeter group
$$W=I_2(m)=\langle s_1,s_2\ |\ s_i^2=1,\
\underbrace{s_1s_2s_1\cdots}_{m\text{ terms }}=
\underbrace{s_2s_1s_2\cdots}_{m\text{ terms }}\rangle.$$
We have
	$$C_W=\langle \gamma_{s_1},\gamma_{s_2},\gamma_S\ |\ 
	\gamma_{s_i}^2=\gamma_S^2=1,\
	\gamma_S\gamma_{s_i}\gamma_S=
\begin{cases}
	\gamma_{3-i} & \text{ if }m \text{ is odd} \\
	\gamma_i & \text{ otherwise}
\end{cases}
\ \rangle$$
and
	\begin{align*}
		f(\gamma_{s_i})&=-t_1+t_{s_i}-t_{s_{3-i}}+
		\frac{2v}{1+v^2}t_{s_is_{3-i}}+t_{w_0}\\
		f(\gamma_S)&=(-1)^mt_1-t_{w_0s_1}-t_{w_0s_2}+t_{w_0}.
	\end{align*}
\end{example}


\begin{thebibliography}{BrMaRouab}
\bibitem[Bo1]{Bo1} C.~Bonnaf\'e,
	``Kazhdan-Lusztig Cells with Unequal Parameters'',
		Springer, 2017.
\bibitem[Bo2]{Bo2} C.~Bonnaf\'e,
	{\em Cells and Cacti},
	Int. Math. Res. Not. {\bf 19} (2016), 5775--5800.
\bibitem[BoRou1]{BoRou1} C.~Bonnaf\'e and R.~Rouquier,
	{\em Calogero-Moser versus Kazhdan-Lusztig cells},
	Pacific J. of Math. {\bf 261} (2013), 45--51.
\bibitem[BoRou2]{BoRou2} C.~Bonnaf\'e and R.~Rouquier,
	{\em Cherednik algebras and Calogero-Moser cells},
	preprint arXiv:1708.09764.
\bibitem[DaJaSc]{DaJaSc} M.~Davis, T.~Januszkiewicz and R.~Scott,
	{\em Fundamental groups of blow-ups},
	Adv. Math. {\bf 177} (2003), 115--179.
\bibitem[De]{De} S.~Devadoss,
	{\em Tessellations of moduli spaces and the mosaic operad},
	in Homotopy Invarient Algebraic Structures (Baltimore, 1998),
		Contemp. Math. 239, Amer.  Math. Soc. 1999, 91--114.
\bibitem[Dr]{Dr} V.~G.~Drinfeld,
	{\em Quasi-Hopf algebras},
		Leningrad Math. J. {\bf 1} (1990), 1419-–1457.
\bibitem[HeKa]{HeKa} A.~Henriques and  J.~Kamnitzer,
	{\em Crystals and coboundary categories},
		Duke Math. J. {\bf 132} (2006), 191--216.
\bibitem[EtHeKaRa]{EtHeKaRa} P.~Etingof, A.~Henriques, J.~Kamnitzer and
	E.~M.~Rains,
	{\em The cohomology ring of the real locus of the moduli
	space of stable curves of genus $0$ with marked points},
	Duke Math. J. {\bf 132} (2006), 191--216.
\bibitem[Lo]{Lo} I.~Losev,
	{\em Cacti and cells}, J. Eur. Math. Soc. {\bf 21} (2019), 
		1729--1750.
\bibitem[Lu1]{Lu1} G.~Lusztig,
	{\em On a theorem of Benson and Curtis},
	J. of Alg. {\bf 71} (1981), 490--498.
\bibitem[Lu2]{Lu2} G.~Lusztig,
	{\em Action of longest element on a Hecke algebra cell module},
	Pacicic Journal of Math. {\bf 279} (2015), 383--396.
\bibitem[Ma]{Ma} A.~Mathas,
	{\em On the left cell representations of Iwahori-Hecke algebras
	of finite Coxeter groups},
	J. London Math. Soc. {\bf 54} (1996), 475--488.
\bibitem[Yu]{Yu} R.~Yu,
	{\em Linearity of Generalized Cactus Groups},
		J. of Alg. {\bf 635} (2023), 256--270.
\end{thebibliography}
\end{document}